\documentclass[letterpaper, 10 pt, conference]{ieeeconf}
\IEEEoverridecommandlockouts

\usepackage{times}

\usepackage{graphicx}
\usepackage{amsmath}
\usepackage{amsfonts}

\usepackage{longtable,tabularx}
\usepackage{mathtools}
\DeclarePairedDelimiter{\norm}{\lVert}{\rVert}

\usepackage{amssymb}

\usepackage{multicol}

\usepackage{algorithm}
\usepackage{algpseudocode}

\usepackage{microtype}

\usepackage{cite}

\usepackage{xcolor}
\usepackage{xfrac}

\usepackage{hhline}

\usepackage{nicefrac}

\usepackage{booktabs} 

\usepackage{multirow}
\usepackage{subcaption}
\usepackage{caption}

\usepackage{threeparttable}
\usepackage{siunitx}

\begin{document}
\title{Adversarial Reinforcement Learning for Robust Control of Fixed-Wing Aircraft under Model Uncertainty}

\author{Dennis~J.~Marquis, Blake Wilhelm, Devaprakash Muniraj, and
 Mazen~Farhood
\thanks{D. Marquis, B. Wilhelm, and M. Farhood are with the Kevin T. Crofton Department of Aerospace and Ocean Engineering, Virginia Tech, Blacksburg, VA 24061, USA (e-mail: \{dennisjm, bawilhelm, farhood\}@vt.edu).}
\thanks{D. Muniraj is with the Department of Aerospace Engineering, Indian Institute of Technology Madras, Chennai 600036, India (email: deva@iitm.ac.in).}}

\markboth{Journal of \LaTeX\ Class Files,~Vol.~14, No.~8, August~2015}%
{Shell \MakeLowercase{\textit{et al.}}: Bare Demo of IEEEtran.cls for IEEE Journals}

\maketitle

\begin{abstract}
This paper presents a reinforcement learning-based path-following controller for a fixed-wing small uncrewed aircraft system (sUAS) that is robust to uncertainties in the aerodynamic model of the sUAS. The controller is trained using the Robust Adversarial Reinforcement Learning framework, where an adversary perturbs the environment (aerodynamic model) to expose the agent (sUAS) to demanding scenarios. In our formulation, the adversary introduces rate-bounded perturbations to the aerodynamic model coefficients. We demonstrate that adversarial training improves robustness compared to controllers trained using stochastic model uncertainty. The learned controller is also benchmarked against a switched uncertain initial condition controller. The effectiveness of the approach is validated through high-fidelity simulations using a realistic six-degree-of-freedom fixed-wing aircraft model, showing accurate and robust path-following performance under a variety of uncertain aerodynamic conditions.
\end{abstract}

\IEEEpeerreviewmaketitle

\section{Introduction}
\label{sec:intro}
Reinforcement learning (RL) is an emerging tool for designing control policies for complex systems. By learning through trial-and-error interactions with an environment, RL can discover strategies that are difficult to derive analytically, particularly for systems with nonlinear dynamics, time-varying disturbances, saturation, or uncertain parameters. For small uncrewed aircraft systems (sUAS), RL has been applied to tasks ranging from low-level attitude stabilization and velocity tracking to high-level planning\cite{Abbeel2006,Hwangbo2017,Koch2019,Kaufmann2022,Dhuheir2023}. Recent surveys of the use of RL for sUAS tasks~\cite{Azar2021,Kurunathan2024} highlight the rapid growth of this field and summarize challenges such as sample efficiency, stability, and sensitivity to discrepancies between simulation and reality, the `sim-to-real gap.'

Applying RL to real-world sUAS remains a challenge due to low sample efficiency and the need for large amounts of training data. A key factor affecting performance is model uncertainty, which can arise from several sources for sUAS: environmental disturbances such as wind, errors in system identification including incorrect modeling of inertial parameters, errors due to linearization for linear controllers, or unmodeled nonlinear effects when operating outside the nominal flight envelope. Prior work demonstrates that, in the presence of these uncertainties, RL can achieve strong performance compared to standard baselines such as PID controllers for fixed-wing sUAS attitude stabilization~\cite{Bohn2019, Bohn2024} and velocity tracking~\cite{Xu2019}. However, maintaining a desired trajectory or path requires the sUAS to counteract environmental disturbances, such as wind, which makes trajectory-tracking and path-following inherently more challenging than lower-level attitude or velocity control~\cite{Muniraj2017,Yang2021}.

Several strategies have been proposed to improve the efficiency and robustness of RL. Curriculum learning~\cite{Bengio2009} gradually increases task difficulty, allowing the agent to learn basic skills before facing more complex scenarios. Domain randomization~\cite{Tobin2017} exposes the policy to a wide range of simulated disturbances and model variations, helping to bridge the sim-to-real gap. In parallel, recent theoretical work integrates RL with control-theoretic principles, providing stronger guarantees on policy performance and convergence~\cite{Hu2023}.

Adversarial approaches have gained traction as a means of enhancing robustness in control and learning. For example, adversarial environments have been leveraged to help learn optimal multi-agent control policies~\cite{Pattanaik2017,Huang2017,Muniraj2018,Lin2019}. The Robust Adversarial Reinforcement Learning (RARL) framework~\cite{Pinto2017} implements this idea by enabling the adversary to have its own learnable policy, which actively perturbs the system to create challenging scenarios for the agent. This formulation directs learning toward difficult scenarios and has been shown to produce policies that generalize better to unforeseen disturbances than those trained with stochastic disturbances on a series of benchmarking OpenAI environments. While RARL has been applied to sUAS planning and high-level trajectory tasks in the past~\cite{Wang2024}, there remains a gap in understanding its effectiveness when it comes to trajectory-tracking and path-following for fixed-wing sUAS under unforseen aerodynamic uncertainties and realistic nonlinear dynamics.

We employ the RARL framework to train an RL-based path-following controller for a fixed wing sUAS against an adversary that introduces rate-bounded perturbations to the coefficients of the aerodynamic model. We evaluate the effectiveness of this approach by comparing policies obtained from adversarial training to those obtained under stochastic aerodynamic uncertainties. In addition, we compare the RL-based controller performance against a switched uncertain initial condition (UIC) controller synthesized using techniques developed in~\cite{Farhood2008}, serving as a model-based benchmark. Simulation results using our realistic sUAS model demonstrate that adversarial training enhances learning performance and yields controllers capable of accurate and robust path-following performance under a variety of model~uncertainties.

The contributions of this work are as follows. First, we extend RL to trajectory-tracking and path-following for fixed-wing sUAS using a realistic six-degree-of-freedom (6DOF) model, rather than simplified kinematic or planning models often used in prior work. Second, to the best of our knowledge, this is the first application of RARL for training stabilizing controllers for sUAS, demonstrating that adversarially trained controllers can generalize better to unforeseen aerodynamic disturbances compared to controllers trained with stochastic perturbations. Third, we present concrete examples of action and observation space selection for the controller and adversary, as well as reward shaping strategies that lead to effective controller performance in the context of path-following and trajectory-tracking.

This paper is organized as follows. Section~\ref{sec:prelim} introduces key preliminaries, including the policy optimization method used in conjunction with RARL and the dynamic model for the fixed-wing sUAS. Section~\ref{sec:training} describes the methodology, including the simulation environment, reference path generation, action and observation spaces, reward shaping, and training procedure. Section~\ref{sec:evaluation} presents an evaluation of the trained controllers, first comparing the performance of controllers trained with adversarial versus stochastic model uncertainty, and then benchmarking the adversarially trained controller against our baseline controller. Finally, Section~\ref{sec:conclude} concludes the paper.

\section{Preliminaries}
\label{sec:prelim}

\subsection{Proximal Policy Optimization}
In RL, an agent interacts with an environment to learn a policy that maximizes a problem-specific cumulative reward. At each discrete time step $k$, the agent makes an observation $\mathbf{o}_k$, selects an action $\mathbf{a}_k$, and obtains a reward $R_k$. For complex systems, $\mathbf{o}_k$ typically includes measurements and known parameters describing the current system and environment configuration. Additionally, a practical optimization algorithm is required to update the policy. In this work, we use Proximal Policy Optimization (PPO)~\cite{schulman2017}, a widely adopted policy gradient method that iteratively updates a policy $\pi$ with parameters $\vartheta$ to maximize the expected reward while avoiding large updates. Inspired by Trust Region Policy Optimization (TRPO)\cite{Schulman2015}, PPO uses a first-order gradient method with the following clipped surrogate objective:
\begin{equation*}
L^{CLIP} = \mathbb{E}\Big[\min\big(\rho_k(\vartheta) \hat{A}_k, \mathrm{clip}(\rho_k(\vartheta), 1-\epsilon, 1+\epsilon) \hat{A}_k\big)\Big],
\end{equation*}
where $\rho_k(\vartheta) = \pi_\vartheta(\mathbf{a}_k|\mathbf{o}_k)/\pi_{\vartheta_\text{old}}(\mathbf{a}_k|\mathbf{o}_k)$, $\hat{A}_k$ is the advantage estimate, and $\epsilon$ is a small clipping parameter. The expectation operator $\mathbb{E}[\cdot]$ denotes the average over trajectories collected from the environment. The function $\mathrm{clip}(x,a,b)$ limits $x$ to the interval $[a,b]$, i.e., $\mathrm{clip}(x,a,b) = \min(\max(x,a),b)$.

A rollout consists of a sequence of consecutive tuples $(\mathbf{o}_k, \mathbf{a}_k, R_k)$ collected by executing the current policy in the environment. The policy parameters $\vartheta$ are then updated by performing gradient ascent on $L^{CLIP}$ using mini-batches of these rollouts. In practice, additional terms for value function error and policy entropy are included to improve stability and encourage exploration.

\subsection{Robust Adversarial Reinforcement Learning} 
While PPO alone can achieve high performance in sUAS applications \cite{Koch2019, Bohn2019}, training against stochastic disturbances often requires large amounts of data to achieve robustness. RARL addresses this problem by optimizing the policy against a trainable adversary, treating uncertainties as adversarial entities. This process produces policies that are not only robust but also generalize better to unforeseen disturbances outside the training distribution. RARL models a two-agent zero-sum game, where a protagonist learns a resilient control policy $\mu$ with parameters $\vartheta^\mu$, while an adversary learns a destabilizing policy $\eta$ with parameters $\vartheta^\eta$.

Training alternates: the adversary's policy is held fixed while the protagonist's policy is updated using an optimizer (here, PPO), then roles are reversed. At instant $k$, the joint environment interaction is $(\mathbf{o}_k^\mu, \mathbf{o}_k^\eta, \mathbf{a}_k^\mu, \mathbf{a}_k^\eta, R_k^\mu, R_k^\eta)$, with $R_k^\mu = -\,R_k^\eta$ for the zero-sum formulation.

\subsection{sUAS Dynamics}
\label{sec:dynamics}
This subsection describes the 6DOF nonlinear model of the CZ-150 sUAS, which is built on the E-flite Carbon-Z Cessna 150T platform. All geometric, inertial, aerodynamic, and actuator parameters were obtained experimentally in~\cite{Marquis2025} and are provided in Appendix~\ref{app:model} for completeness. 

The gravitational constant is denoted by $g$. The vehicle mass, mean aerodynamic chord, wingspan, and wing reference area are represented by $m$, $\bar{c}$, $b$, and $S$, respectively. Assuming symmetry in the $x$–$z$ plane, the inertia tensor only contains the nonzero terms $J_{xx}, J_{yy}, J_{zz}, J_{xz}$, while $J_{xy}$ and $J_{yz}$ vanish.

The position of the vehicle in the inertial North–East–Down frame $\mathcal{F}_I$ is $\mathbf{p} = [x \; y \; z]^T$. Linear and angular velocities, both expressed in the body-fixed frame $\mathcal{F}_b$, are $\mathbf{v} = [u \; v \; w]^T$ and $\boldsymbol{\omega} = [p \; q \; r]^T$, respectively. The attitude vector $\boldsymbol{\theta} = [\phi \; \theta \; \psi]^T$ parameterizes the rotation matrix $R_{Ib}(\boldsymbol{\theta})$, from $\mathcal{F}_b$ to $\mathcal{F}_I$, and the Euler rate matrix $\varepsilon(\phi,\theta)$:
\begin{equation*}
R_{Ib}(\boldsymbol{\theta}) \coloneqq
\begin{bmatrix}
c_\theta c_\psi & s_\phi s_\theta c_\psi - c_\phi s_\psi & c_\phi s_\theta c_\psi + s_\phi s_\psi \\
c_\theta s_\psi & s_\phi s_\theta s_\psi + c_\phi c_\psi & c_\phi s_\theta s_\psi - s_\phi c_\psi \\
-s_\theta & s_\phi c_\theta & c_\phi c_\theta
\end{bmatrix},
\end{equation*}
\begin{equation*}
\varepsilon(\phi,\theta) \coloneqq
\begin{bmatrix}
1 & s_\phi t_\theta & c_\phi t_\theta \\
0 & c_\phi & -s_\phi \\
0 & s_\phi/c_\theta & c_\phi/c_\theta
\end{bmatrix},
\end{equation*}
where $s_x \coloneqq \sin{x}$, $c_x \coloneqq \cos{x}$, and $t_x \coloneqq \tan{x}$.

The control vector is $\boldsymbol{\delta} = [\delta_E \; \delta_A \; \delta_R \; \delta_T]^T$, corresponding to elevator, aileron, rudder, and throttle commands. Surface deflections are measured in radians, while the throttle input is expressed through the motor’s rotation rate in revolutions per second. The inertial-frame wind velocity is $\mathbf{v}_w = [v_{w,x} \; v_{w,y} \; v_{w,z}]^T$, and the relative velocity of the aircraft with respect to the air is $\mathbf{v}_r = [u_r \; v_r \; w_r]^T = \mathbf{v} - R_{Ib}^T(\boldsymbol{\theta})\mathbf{v}_w$. The forces and moments, $\mathbf{F}(\mathbf{v}_r,\boldsymbol{\omega},\boldsymbol{\delta})$ and $\mathbf{M}(\mathbf{v}_r,\boldsymbol{\omega},\boldsymbol{\delta})$, are here expressed simply as $\mathbf{F}(t)$ and $\mathbf{M}(t)$. The equations of motion for the 6DOF system are then given by
\begin{equation*}
\begin{aligned}
\dot{\mathbf{p}} &= R_{Ib}(\boldsymbol{\theta})\mathbf{v}, \quad
\dot{\boldsymbol{\theta}} = \varepsilon(\phi,\theta)\boldsymbol{\omega}, \\
\dot{\mathbf{v}} &= \mathbf{v}\times \boldsymbol{\omega} + R_{Ib}^T(\boldsymbol{\theta})[0 \; 0 \; g]^T + \tfrac{1}{m}\mathbf{F}(t), \\
\dot{\boldsymbol{\omega}} &= J^{-1}\left[J\boldsymbol{\omega}\times \boldsymbol{\omega} + \mathbf{M}(t)\right].
\end{aligned}
\end{equation*}

Forces and moments are parameterized through the aerodynamic coefficients as
\begin{equation*}
\begin{aligned}
\mathbf{F}(t) &=
\begin{bmatrix}
C_X(t)\,\bar{q}S \\
C_Y(t)\,\bar{q}S \\
C_Z(t)\,\bar{q}S
\end{bmatrix}, \quad
\mathbf{M}(t) =
\begin{bmatrix}
C_L(t)\,\bar{q}Sb \\
C_M(t)\,\bar{q}S\bar{c} \\
C_N(t)\,\bar{q}Sb
\end{bmatrix}.
\end{aligned}
\end{equation*}

The dynamic pressure is $\bar{q} = \tfrac{1}{2}\rho_\mathrm{air}V_r^2$, $\rho_\mathrm{air}$ is the air density and $V_r = \norm{\mathbf{v}_r}$ is the airspeed, where $\norm{\cdot}$ denotes the Euclidean norm. The regressor variables in the aerodynamic model include the angle of attack $\alpha = \arctan(w_r/u_r)$, sideslip angle $\beta = \arcsin(v_r/V_r)$, nondimensional angular rates $\hat{p} = pb/(2V_r)$, $\hat{q} = q\bar{c}/(2V_r)$, $\hat{r} = rb/(2V_r)$, control inputs $\delta_E$, $\delta_A$, $\delta_R$, and the inverse advance ratio $\mathcal{J} = \delta_T D_\mathrm{prop}/V_r$, with $D_\mathrm{prop}$ as the propeller diameter.

The commanded inputs $\boldsymbol{\delta}^\mathrm{cmd} = [\delta_E^\mathrm{cmd} \; \delta_A^\mathrm{cmd} \; \delta_R^\mathrm{cmd} \; \delta_T^\mathrm{cmd}]^T$ are given in radians for control surfaces and in revolutions per second for throttle. These values are obtained from centered PWM signals $\boldsymbol{\delta}^\mathrm{PWM}$ using static polynomial mappings (cubic for control surfaces, linear for throttle). The final actuator states $\boldsymbol{\delta}$ are then produced by passing $\boldsymbol{\delta}^\mathrm{cmd}$ through first-order actuator dynamics.

To enable path-following, we adopt the virtual vehicle framework introduced in \cite{Fry2020} (see Sec. II-B, Eqs. (6–11)), in which the sUAS dynamics are expressed relative to a virtual vehicle moving along the desired path and controlled based on the error with respect to this reference.

\section{Training Methodology}
\label{sec:training}
This section describes the framework used to train an RL-based path-following controller with RARL, including the simulation setup, reference path generation, action and observation spaces, and reward design.

\subsection{Simulation Environment}
\label{sec:sim}
To enable training, a custom Python environment has been developed based on the CZ-150 vehicle dynamics. To make the simulated behavior more representative of actual flight, realistic disturbances are incorporated into the simulation. Gaussian measurement noise is applied to all sensor channels, with standard deviations of $0.01$~rad/s for $\boldsymbol{\omega}$, $2$~m/s for $V_r$, $0.01$~rad for $\phi$ and $\theta$, $0.1$~rad for $\psi$, $0.03$~m for $x$ and $y$, $0.01$~m for $z$, and $0.03$~m/s$^2$ for translational acceleration in the body frame $\mathbf{f}_b = [a_x \; a_y \;a_z]^T$. Wind is modeled as a steady wind component with a randomly sampled magnitude in $[3$ m/s$, 7$ m/s$]$ and a heading in $[0^\circ,360^\circ]$ plus wind gusts generated using the Dryden turbulence model~\cite{Real1993} with a moderate reference wind speed of 30~knots at low altitude. Controller inputs are delayed up to one simulation step. The continuous-time dynamics described in Section~\ref{sec:dynamics} are integrated using a fourth-order Runge--Kutta (RK4) scheme with a fixed step size of $\Delta t = 0.04$~s. Since the implementation is in discrete time with the sampling period $\Delta t$, we write 
\begin{equation*}
t_k = k \Delta t, \quad k = 0,1,\dots,N,
\end{equation*}
where $N = t_f/\Delta t$. For any continuous-time variable $x(t)$, its discrete-time sequence is defined as $x_k = x(t_k)$.

In addition to stochastic disturbances, an adversary is introduced in the form of aerodynamic model uncertainty. Each aerodynamic coefficient $C_{i,k}$, for $i \in \{X,Y,Z,L,M,N\}$, is perturbed as $C_{i,k} \rightarrow C_{i,k} + \Delta_{C_i,k}$, defining the perturbation vector $\boldsymbol{\Delta}_{\mathbf{C},k} \in \mathbb{R}^6$.
The perturbations are constrained by
\begin{equation*}
\boldsymbol{\Delta}_{\mathbf{C}}^- \leq \boldsymbol{\Delta}_{\mathbf{C},k} \leq \boldsymbol{\Delta}_{\mathbf{C}}^+,\quad 
\boldsymbol{\nu}_{\mathbf{C}}^- \leq \boldsymbol{\Delta}_{\mathbf{C},k} - \boldsymbol{\Delta}_{\mathbf{C},k-1} \leq \boldsymbol{\nu}_{\mathbf{C}}^+,
\end{equation*}
where the bounds $\boldsymbol{\Delta}_{\mathbf{C}}^-$, $\boldsymbol{\Delta}_{\mathbf{C}}^+$, $\boldsymbol{\nu}_{\mathbf{C}}^-$, and $\boldsymbol{\nu}_{\mathbf{C}}^+$ are determined from flight data. More details regarding the adversary, including its implementation and training, are provided in the sequel.

\subsection{Reference Path Generation}
\label{sec:ref}
Lemniscate-like reference paths are generated by concatenating motion primitives, using an approach presented in~\cite{Arifianto2015}. Each path consists of a straight-and-level segment, a coordinated ascending or descending turn parameterized by the inverse radius of curvature $\kappa$ and flight path angle $\gamma$, another straight-and-level segment, and a corresponding coordinated descending or ascending turn to complete the loop. Each segment corresponds to a distinct trim, defined by the equilibrium values of the states $\phi_\mathrm{trim}$, $\theta_\mathrm{trim}$, $\mathbf{v}_\mathrm{trim}$, $\boldsymbol{\omega}_\mathrm{trim}$, $\boldsymbol{\delta}_\mathrm{trim}$, and the input $\boldsymbol{\delta}_\mathrm{trim}^\mathrm{cmd}$, obtained at the nominal airspeed $V_\mathrm{r,nom} = 21$~m/s using a nonlinear solver. For each segment, the trims are used to generate time-varying reference histories of $x$, $y$, $z$, and $\psi$ by integrating the vehicle dynamics. These segment references are then concatenated, yielding a reference $\mathbf{x}_\mathrm{ref}$ and corresponding input $\boldsymbol{\delta}_\mathrm{ref}^\mathrm{cmd}$ for each lemniscate path, where the full state is $\mathbf{x} = [\mathbf{p}^T \; \mathbf{v}^T \; \boldsymbol{\theta}^T \; \boldsymbol{\omega}^T \; \boldsymbol{\delta}^T]^T$.

The coordinated turns are based on a predefined set of trims, each defined by a pair of $(\kappa, \gamma) \in K_\mathrm{ref} \times \Gamma_\mathrm{ref}$, where
\begin{equation*}
\begin{gathered}
K_\mathrm{ref} = \{-0.02, -0.012, 0.012, 0.02\}, \\
\Gamma_\mathrm{ref} = \{-0.21, -0.11, 0, 0.11, 0.21\},
\end{gathered}
\end{equation*}
while $\kappa = \gamma = 0$ is used to generate straight-and-level flight segments. This yields 20 distinct coordinated-turn patterns, each combined with the straight segments to produce a respective lemniscate-like path.

Although in this work we focus on evaluating a path-following controller, the architecture is compatible with trajectory-tracking. In high-wind conditions, predefined trajectories may become dynamically infeasible. The path-following approach relaxes the timing constraints, though wind disturbances still make accurate path-following challenging.

\subsection{Action and Observation Spaces}
In this work, the policies of the protagonist and adversary are each represented by a neural network with two hidden layers of 64 neurons and Tanh activations. The protagonist's action space consists of control perturbations relative to the reference input. The adversary's action space consists of perturbations to the aerodynamic coefficients that evolve according to the previously defined rate bounds.  
The respective actions are given by
\begin{equation*}
\mathbf{a}_k^\mu = \boldsymbol{\delta}_k^\mathrm{cmd} - \boldsymbol{\delta}_\mathrm{ref,k}^\mathrm{cmd} \ \textrm{and} 
\ 
\mathbf{a}_k^\eta = \boldsymbol{\Delta}_{\mathbf{C},k} - \boldsymbol{\Delta}_{\mathbf{C},k-1}.
\end{equation*}

The measurement vector is $\mathbf{y} = [\boldsymbol{\omega}^T \; V_r \; \boldsymbol{\theta}^T \; \mathbf{p}^T \; \mathbf{f}_b^T]^T$. The protagonist's  observation vector $\mathbf{o}_k^\mu$ consists of the tracking error $\bar{\mathbf{y}}_k = \mathbf{y}_k - \mathbf{y}_{\mathrm{ref},k}$, the reference input $\boldsymbol{\delta}_\mathrm{ref,k}^\mathrm{cmd}$, the previous control input $\boldsymbol{\delta}_{k-1}^\mathrm{cmd}$, the control margin $\mathbf{m}_k$, and the maneuver scheduling parameters $\kappa_k$ and $\gamma_k$. The control margin is
\begin{equation*}
\mathbf{m}_k = \min \Bigg\{
\max\Big(\mathbf{0}, \tfrac{\boldsymbol{\delta}_{\max} - \boldsymbol{\delta}_k^\mathrm{cmd}}{\boldsymbol{\delta}_{\max} - \boldsymbol{\delta}^\mathrm{cmd}_\mathrm{ref,k}}\Big),
\max\Big(\mathbf{0}, \tfrac{\boldsymbol{\delta}_k^\mathrm{cmd} - \boldsymbol{\delta}_{\min}}{\boldsymbol{\delta}^\mathrm{cmd}_\mathrm{ref,k} - \boldsymbol{\delta}_{\min}}\Big)
\Bigg\},
\end{equation*}
where $\boldsymbol{\delta}_{\max}$ and $\boldsymbol{\delta}_{\min}$ are the respective upper and lower actuator saturation limits and the max and min functions are defined component-wise. By construction, each element of $\mathbf{m}_k$ is 1 at trim and 0 at saturation. This observation space provides the agent with the tracking error to correct deviations from the reference, context for temporal control decisions via previous inputs, actuator saturation limits through the control margin, and scheduling parameters allowing the controller to generalize across the family of reference trajectories.

The adversary's observation vector $\mathbf{o}_k^\eta$ consists of the same tracking error $\bar{\mathbf{y}}_k$ and the current aerodynamic perturbation $\Delta_{\mathbf{C},k}$, allowing the adversary to efficiently learn its effective saturation bounds. All action and observation spaces are normalized to $[-1,1]$ for improved neural network performance, although the physical ranges are presented here for clarity.

\subsection{Reward Shaping}
The total reward at each instant $k$ is
\begin{equation*}
R_k^\mu = -R_k^\eta = R_{\mathrm{tracking},k} + R_{\mathrm{input},k}.
\end{equation*}
The protagonist is rewarded for maintaining accurate position, attitude, and body rates, using normalized exponential functions of the tracking errors, given by
\begin{equation*}
R_{\mathrm{tracking},k} = \sum_{i=1}^{9} k_{1,i} \, e^{-k_{2,i} |\bar{y}_{k,i}|},
\end{equation*}
\noindent
where $i=1,2,3$ correspond to body rates $p,q,r$ with $k_{1,i}=0.1$, $k_{2,i}=1.0$; 
$i=4,5,6$ correspond to attitudes $\phi,\theta,\psi$ with $k_{1,i}=0.2$, $k_{2,i}=5.0$; 
and $i=7,8,9$ correspond to positions $x,y,z$ with $k_{1,i}=0.5$, $k_{2,i}=0.37$. 
Here, $\bar{y}_{k,i}$ denotes the tracking error of output $i$ at instant $k$. The tracking rewards are shaped to incentivize keeping errors below particular bounds and to avoid punishing small deviations. The input reward combines a barrier term that discourages control saturation and a rate term that discourages rapid input changes:
\begin{equation*}
R_{\mathrm{input},k} 
= k_3 \sum_{j=1}^{4} \log\!\big(m_{k,j} + 10^{-6}\big)
- k_4 \, \norm{\boldsymbol{\delta}_k^\mathrm{cmd} - \boldsymbol{\delta}_{k-1}^\mathrm{cmd}}^2,
\end{equation*}
\noindent
where $k_3=0.05$, $k_4=0.2$, $j=1,2,3,4$ correspond to $\delta_E^\mathrm{cmd}, \delta_A^\mathrm{cmd}, \delta_R^\mathrm{cmd}, \delta_T^\mathrm{cmd}$, respectively, and $m_{k,j}$ is the control margin for input $j$ at instant 
$k$. Policies obtained using RL often converge to high-gain bang-bang control schemes, which can cause excessive actuator wear and perform worse in practice than in simulation; including a rate term helps smooth the resulting control signals.

\subsection{Training}
Training of all RARL-based path-following controllers is performed on Virginia Tech's Advanced Research Computing cluster. We implement PPO using Stable-Baselines3~\cite{Raffin2021StableBaselines3} with PyTorch~\cite{paszke2019pytorch}. Each policy update is based on rollouts collected from $n_\mathrm{envs}$ parallel environments, increasing the batch size per iteration and reducing variance in the gradient estimates. The RARL procedure is summarized in Algorithm~\ref{alg:RARL}. At each iteration, the protagonist is updated first while the adversary is held fixed, followed by an update of the adversary with the protagonist fixed. Each update uses rollouts collected from all parallel environments, yielding $N_\mathrm{roll} = n_\mathrm{envs} \times n_\mathrm{steps}$ steps per iteration, over $N_\mathrm{iter}$ iterations. In Algorithm~\ref{alg:RARL}, the function $\mathrm{PPO}(\text{rollouts}, \pi, \vartheta)$ performs a standard PPO update, adjusting the parameters $\vartheta$ of policy $\pi$ via gradient ascent, using the collected rollouts.

Key hyperparameters used for training the PPO models are summarized in Table~\ref{table:training_params}. These values are commonly used in continuous-control tasks and have been shown to produce stable learning in similar RL settings.

\begin{table}[h!]
\centering
\caption{Training Hyperparameters}
\label{table:training_params}
\begin{tabular}{@{}l|l|l|l@{}}
\toprule
\textbf{Hyperparameter} & \textbf{Value} & \textbf{Hyperparameter} & \textbf{Value} \\ \midrule
$N_\mathrm{iter}$ & \{100, 500, 1000\} & $n_\mathrm{envs}$ & 8 \\
learning rate & $3\times10^{-4}$ & $n_\mathrm{steps}$ & 2048 \\
batch size & 64 & $n_\mathrm{epochs}$ & 10 \\
discount factor & 0.99 & GAE parameter & 0.95 \\
$\epsilon$ (clip range) & 0.2 & entropy coefficient & 0 \\
value loss coefficient & 0.5 & gradient clipping & 0.5 \\ \bottomrule
\end{tabular}
\end{table}

\begin{algorithm}
\caption{RARL Procedure~\cite{Pinto2017}}\label{alg:RARL}
\begin{algorithmic}[1]
  \State \textbf{Input:} Stochastic policies $\mu$ and $\eta$
  \State \textbf{Initialize:} Parameters $\vartheta_0^{\mu}$, $\vartheta_0^{\eta}$
  \For{$i = 1,2,...,N_\text{iter}$}
    \For{each policy $\pi \in \{\mu,\eta\}$}
      \State Update $\vartheta_i^\pi$ with the other policy fixed:
      \For{$j = 1,...,N_\mathrm{roll}$}
        \State Collect rollouts $(\mathbf{o}_k^\mu, \mathbf{o}_k^\eta, \mathbf{a}_k^{\mu}, \mathbf{a}_k^{\eta}, R_k^{\mu}, R_k^{\eta})$
        \State $\vartheta_i^\pi \gets \mathrm{PPO}(\text{rollouts}, \pi, \vartheta_i^\pi)$
      \EndFor
    \EndFor
  \EndFor
  \State \textbf{Return:} $\vartheta_{N_\text{iter}}^{\mu}, \vartheta_{N_\text{iter}}^{\eta}$
\end{algorithmic}
\end{algorithm}

\section{Evaluation}
\label{sec:evaluation}
This section evaluates controller performance under uncertain aerodynamic conditions, first comparing stochastic and RARL-trained policies to assess generalization, and then benchmarking a fully trained RL-based controller against the baseline switched UIC controller, whose synthesis and design details are provided in Appendix~\ref{app:uic}.

\subsection{Adversarial versus Stochastic Controller Performance}
To evaluate the RARL framework for training path-following controllers, we train policies under adversarial model uncertainty as described in Section~\ref{sec:training}. For comparison, controllers are also trained under stochastic model uncertainty with the same amplitude and rate bounds. Each type is trained for two horizons, $N_\mathrm{iter} = 100$ and $500$, with four independent policies per configuration (since training itself is stochastic). To ensure fairness in the comparisons, the adversary used in the evaluations of the RL-based controllers is trained against a separate switched UIC controller, thus generating challenging, protagonist-independent disturbances.

\begin{figure*}[t]
\setlength{\abovecaptionskip}{2pt}
\setlength{\belowcaptionskip}{2pt}
    \centering
    \begin{subfigure}[c]{0.32\textwidth}
        \centering
        \includegraphics[scale=0.77]{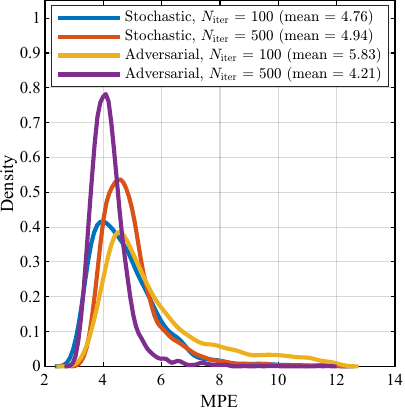}
        \label{fig:compare_mpe}
    \end{subfigure}
    \hfill
    \begin{subfigure}[c]{0.32\textwidth}
        \centering
        \includegraphics[scale=0.77]{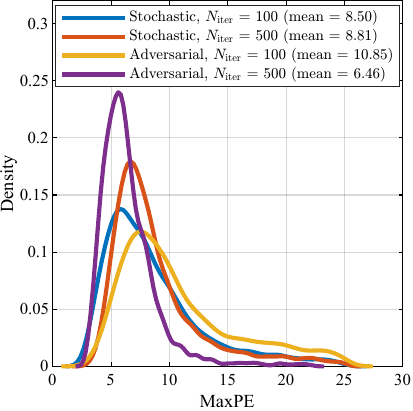}
        \label{fig:compare_maxpe}
    \end{subfigure}
    \hfill
        \begin{subfigure}[c]{0.32\textwidth}
        \centering
        \includegraphics[scale=0.77]{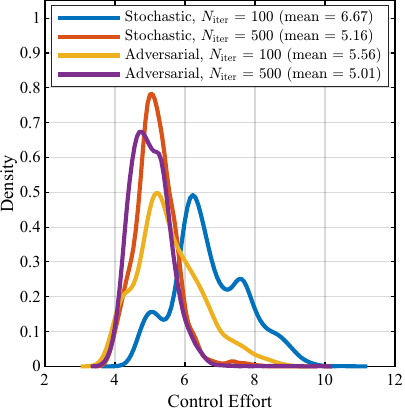}
        \label{fig:compare_ce}
    \end{subfigure}
    \caption{Aggregate of performance metrics for controllers trained under stochastic and adversarial aerodynamic uncertainties.}
    \label{fig:aggregate_metrics}
\end{figure*}

The performance metrics are mean path error (MPE), maximum path error (MaxPE), and control effort, defined as $\sum_{k=1}^{N} ||\mathbf{a}_k^\mu||$, where $N$ is the total number of time steps in the trial. Each policy is evaluated over 1000 trials. Raw trial results are shown in Table~\ref{table:raw_data}, while aggregate results are summarized in Figure~\ref{fig:aggregate_metrics}.

The top-performing controller in all metrics is the adversarially trained policy with $N_\mathrm{iter} = 500$, achieving an average MPE of 4.21\,m, average MaxPE of 6.46\,m, and average control effort of 5.01. Both stochastically trained controllers perform worse than the top-performing adversarially trained controller. The 100-iteration stochastic controller achieves an average MPE of 4.76\,m, average MaxPE of 8.50\,m, and average control effort of 6.67, indicating moderate tracking performance but comparatively high control usage. In contrast, the 500-iteration stochastic controller reaches an average MPE of 4.94\,m, average MaxPE of 8.81\,m, and a lower average control effort of 5.16, showing that extended training reduces the amount of control effort required while slightly increasing the tracking error. This suggests that longer training allows the stochastic policy to find more efficient control strategies, even if tracking precision is marginally compromised.  

\begin{table}[hb!]
\centering
\caption{Raw trial data for trained controllers.$^*$}
\label{table:raw_data}
\begin{tabular}{@{}l|l|l@{}}
\toprule
\textbf{Training Method, $N_\mathrm{iter}$} & \textbf{Metric} & \textbf{Means} \\ \midrule
\multirow{3}{*}{Stochastic, 100} & MPE (m) & (4.58, 4.37, 4.73, 5.34) \\
                                & MaxPE (m) & (8.05, 6.92, 9.99, 9.05) \\
                                & Control Effort & (5.80, 7.40, 6.98, 6.50) \\ \midrule
\multirow{3}{*}{Stochastic, 500} & MPE (m) & (5.31, 4.51, 4.77, 5.15) \\
                                & MaxPE (m) & (10.65, 7.25, 8.17, 9.18) \\
                                & Control Effort & (4.97, 5.24, 4.98, 5.43) \\ \midrule
\multirow{3}{*}{Adversarial, 100} & MPE (m) & (5.65, 5.05, 5.62, 7.01) \\
                                & MaxPE (m) & (10.74, 9.49, 10.18, 13.00) \\
                                & Control Effort & (6.09, 4.77, 5.53, 5.84) \\ \midrule
\multirow{3}{*}{Adversarial, 500} & MPE (m) & (4.45, 4.22, 4.05, 4.12) \\
                                & MaxPE (m) & (6.66, 6.62, 5.78, 6.76) \\
                                & Control Effort & (5.21, 5.35, 4.88, 4.60) \\
\bottomrule
\end{tabular}
\caption*{\footnotesize $^*$Means are computed per policy (1000 trials each).}
\end{table} 

The differences in tracking performance and control effort are closely related to the reward structure. In the current formulation, control rates are relatively heavily penalized, which encourages smoother, lower-amplitude commands at the expense of tighter tracking. Adjusting the reward coefficients could potentially tighten path-following performance, while a sparser reward structure, like the one used in~\cite{Bohn2024}, could also be explored to emphasize achieving the goal rather than strictly penalizing input effort.  

Furthermore, the distributions of MPE and MaxPE in Figure~\ref{fig:aggregate_metrics} reveal that the stochastic controllers exhibit longer tails, indicating there are scenarios where the controller is substantially disrupted. This effect is most pronounced in maximum error: all controllers except the 500-iteration adversarial policy occasionally reach up to 25\,m MaxPE, while the best adversarial controller rarely exceeds 15\,m error. This result illustrates the improved robustness of adversarial training against extreme disturbances.  

The worst-performing controller is the adversarially trained policy with 100 iterations, which likely failed to converge while simultaneously contending with a trainable adversary, achieving 5.83\,m average MPE, 10.85\,m average MaxPE, and 5.56 average control effort. Overall, these results demonstrate that with a sufficient learning horizon (here, 500 iterations), an adversarially trained controller better generalizes to this unforeseen adversarial modeling error and mitigates worst-case deviations compared to its stochastically trained counterparts.

\begin{figure*}[t]
\setlength{\abovecaptionskip}{2pt}
\setlength{\belowcaptionskip}{2pt}
    \centering
    \begin{subfigure}[t]{0.32\textwidth}
        \centering
        \includegraphics[scale=0.77]{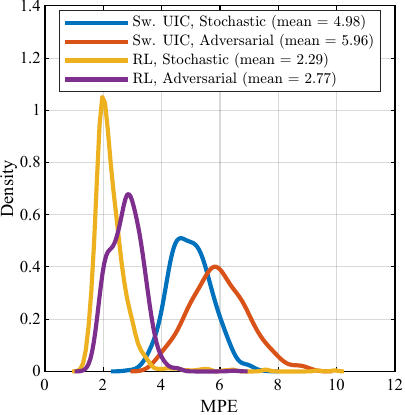}
        \label{fig:baseline_mpe}
    \end{subfigure}
    \hfill
    \begin{subfigure}[t]{0.32\textwidth}
        \centering
        \includegraphics[scale=0.77]{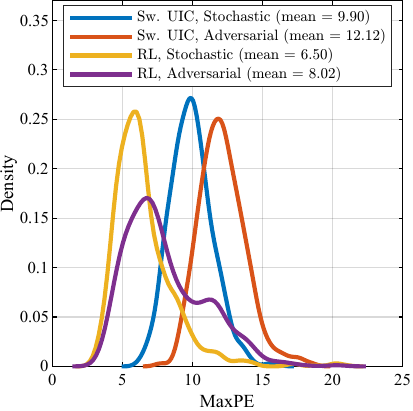}
        \label{fig:baseline_maxpe}
    \end{subfigure}
    \hfill
        \begin{subfigure}[t]{0.32\textwidth}
        \centering
        \includegraphics[scale=0.77]{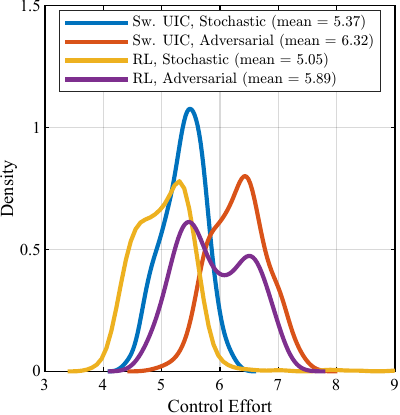}
        \label{fig:baseline_ce}
    \end{subfigure}
    \caption{Switched UIC and RL-based controllers evaluated against stochastic and adversarial aerodynamic uncertainties.}
    \label{fig:baseline}
\end{figure*}

\begin{figure}[!ht]
\centering
\setlength{\abovecaptionskip}{2pt}
\setlength{\belowcaptionskip}{2pt}
\includegraphics[scale=0.81]{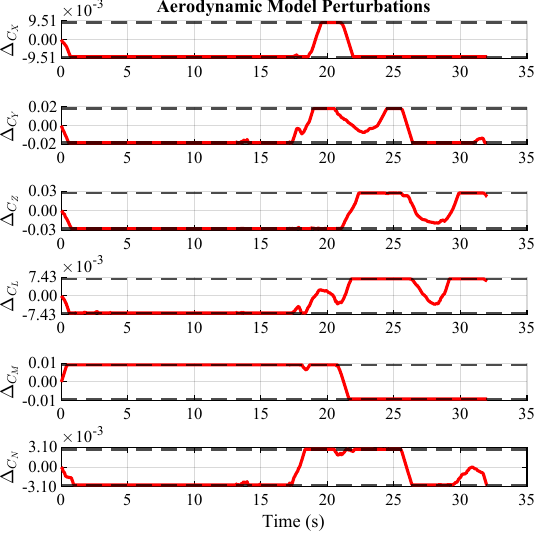}
\caption{Adversarial aerodynamic model perturbations for an example simulation episode.}
\label{fig:sim_deltaC}\vspace{-7mm}
\end{figure}

\subsection{Adversarial versus Baseline Controller Performance}
To benchmark a well-trained RL-based controller, we train a new policy against adversarial model uncertainty for $N_\mathrm{iter}= 1000$ and compare it with the baseline switched UIC controller. Each controller is evaluated under both adversarial and stochastic model uncertainties, with 1000 trials per configuration. The simulations are demanding, with wind speeds reaching up to 7\,m/s.

Figure~\ref{fig:baseline} summarizes the benchmarking results. The RL-based controller objectively outperforms the switched UIC controller across all metrics, despite the highly stressed simulation environment described in Section~\ref{sec:sim}. Averaged over 1000 trials, the RL policy achieves a mean path error of 2.29\,m under stochastic disturbances and 2.77\,m under adversarial disturbances, demonstrating robust tracking even under challenging model perturbations. In comparison, the switched UIC controller exhibits higher errors, with 4.98\,m under stochastic disturbances and 5.96\,m under adversarial disturbances. Maximum path error follows the same trend: the RL controller averages at 6.50\,m (stochastic) and 8.02\,m (adversarial), whereas the switched UIC controller averages at 9.90\,m (stochastic) and 12.12\,m (adversarial), indicating that the RL policy reduces the severity of worst-case deviations. Control effort is lower for the RL controller, with 5.05 (stochastic) and 5.89 (adversarial), compared to 5.37 (stochastic) and 6.32 (adversarial) for the switched UIC controller, reflecting a slightly more efficient use of actuators.

These results also highlight the effectiveness of the adversarial disturbances. Both controllers perform worse under adversarial perturbations than under stochastic ones, confirming that the adversary successfully identifies and exploits system vulnerabilities. Moreover, the RL policy consistently respects actuator limits, rarely exceeding 0.02\% saturation, whereas the switched UIC controller occasionally saturates elevator (up to $\sim$0.6\%) and rudder (up to $\sim$2\%). This improved saturation behavior is largely attributed to the inclusion of the control margin reward term, which encourages the RL-based controller to stay within bounds compared to earlier versions. Overall, the RL-based controller achieves tighter, more efficient, and more robust control, outperforming the switched UIC controller even under substantial disturbances.

\begin{figure*}[t]
\centering
\setlength{\abovecaptionskip}{2pt}
\setlength{\belowcaptionskip}{2pt}
\begin{subfigure}[c]{0.32\textwidth}
    \centering
    \includegraphics[scale=0.8]{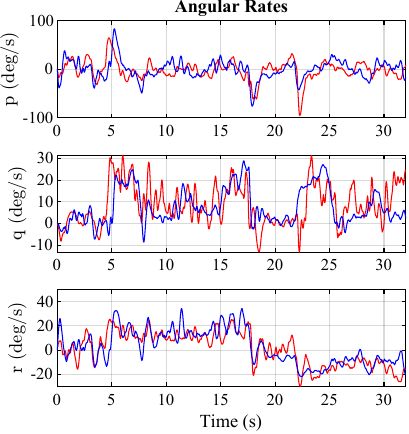}
    \label{fig:sim_rates}
\end{subfigure}
\hfill
\begin{subfigure}[c]{0.32\textwidth}
    \centering
    \includegraphics[scale=0.8]{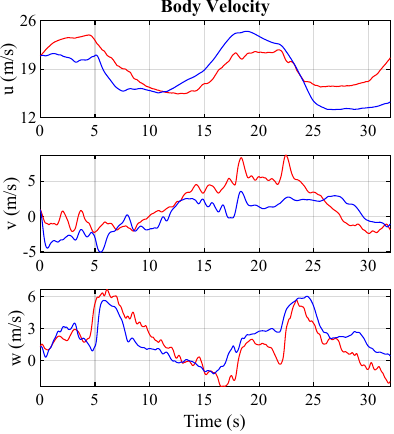}
    \label{fig:sim_vel}
\end{subfigure}
\hfill
\begin{subfigure}[c]{0.32\textwidth}
    \centering
    \includegraphics[scale=0.8]{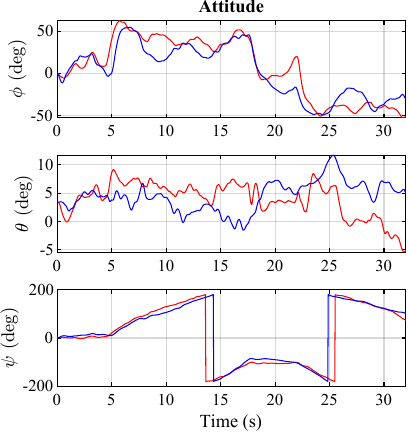}
    \label{fig:sim_attitude}
\end{subfigure}
\vspace{2mm} \\
\begin{subfigure}[c]{0.32\textwidth}
    \centering
    \includegraphics[scale=0.8]{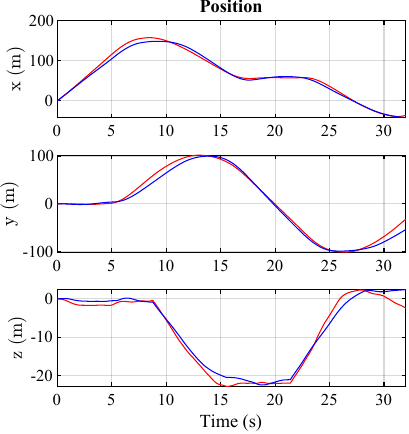}
    \label{fig:sim_pos}
\end{subfigure}
\hfill
\begin{subfigure}[c]{0.32\textwidth}
    \centering
    \includegraphics[scale=0.8]{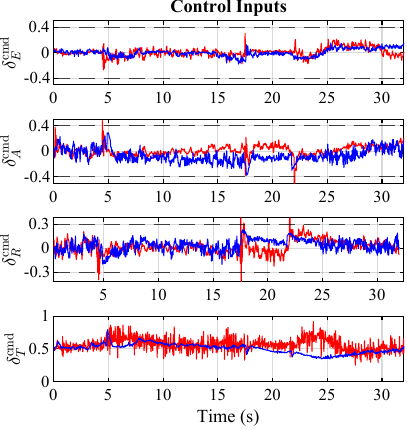}
    \label{fig:sim_control}
\end{subfigure}
\hfill
\begin{subfigure}[c]{0.32\textwidth}
    \centering
    \includegraphics[scale=0.8]{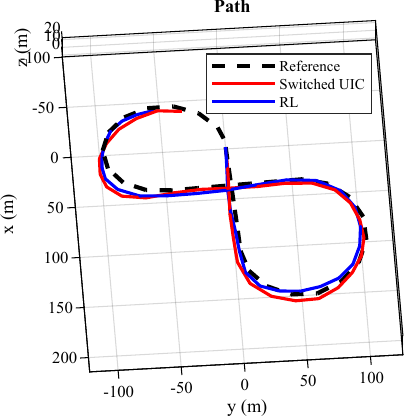}
    \label{fig:sim_path}
\end{subfigure}
\caption{State and control input histories of the switched UIC controller and RL controller for an example simulation episode.}
\label{fig:sim_states}
\end{figure*}

To illustrate controller behavior, we present an example simulation on a lemniscate-like path corresponding to $(\kappa, \gamma) = (0.02, 0.21)$ where adversarial model perturbations are present. The perturbation plot in Figure~\ref{fig:sim_deltaC} shows that the adversarial strategy tends to operate at the perturbation limits, an intuitive outcome since there is no reward disincentive. At the start of the transition to straight and level flight from the first coordinated turn (the midpoint), the adversary adopts a more oscillatory strategy, particularly in $C_Y$ and $C_L$, possibly exploiting the increased vulnerability of the controllers during path transitions and the changes in the sUAS’s direction relative to the wind.

The state and control input plots are presented in Figure~\ref{fig:sim_states}. Angular rate behavior is similar across controllers. The body velocity plots highlight the influence of strong wind on path-following, with $u$ ranging between 12 and 25\,m/s for both controllers. While the lateral velocity $v$ is near zero at trim, it is non-negligible in this simulation, reflecting sideslip induced by wind and control actions. The largest disparity between the two trajectories occurs in the pitch angle $\theta$, where the switched UIC controller produces nose-down attitude at the end of the episode, while the RL controller maintains more positive pitch. The RL controller develops a low-amplitude oscillatory control policy, reflecting its preference for bang-bang actuation constrained by rate limits and the control margin penalty. In contrast, the switched UIC controller demonstrates particularly noisy throttle commands, with instances where both aileron and rudder inputs reach saturation. The paths do not form a complete loop because the simulation executes for a fixed duration corresponding to the nominal time to complete a loop (in the absence of disturbances and perturbations).

Overall, these simulation results demonstrate that the RL-based controller provides consistently improved tracking performance, reduced control effort, and enhanced robustness to both stochastic and adversarial disturbances. The example simulation highlights how the RL policy improves state and control trajectories, illustrating the practical advantages of adversarial training within the path-following framework. Finally, it is important to stress that the obtained results reflect the specific formulations and tuning used in this work and should not be interpreted as general since alternative formulations could produce different outcomes. 

\section{Conclusion}\label{sec:conclude}
This work presents an RL-based path-following controller for a fixed-wing sUAS, trained using the RARL framework to improve robustness to uncertainties. The adversary introduces rate-bounded aerodynamic model perturbations to expose the agent to demanding scenarios. Each controller is evaluated in high-fidelity simulations, and the learned RL-based controller is benchmarked against a switched UIC controller synthesized using the CZ-150 model. Results show that adversarial training can yield superior tracking performance and lower control effort compared to both stochastically trained controllers and the switched UIC baseline. Future work includes applying RARL to additional sources of uncertainty and flight-testing the RL-based controller to verify real-world performance.

\bibliographystyle{ieeetr}

\begin{thebibliography}{10}

\bibitem{Abbeel2006}
P.~Abbeel, A.~Coates, M.~Quigley, and A.~Ng, ``An application of reinforcement learning to aerobatic helicopter flight,'' in {\em Advances in Neural Information Processing Systems}, vol.~19, MIT Press, 2006.

\bibitem{Hwangbo2017}
J.~Hwangbo, I.~Sa, R.~Siegwart, and M.~Hutter, ``Control of a quadrotor with reinforcement learning,'' {\em IEEE Robotics and Automation Letters}, vol.~2, no.~4, pp.~2096--2103, 2017.

\bibitem{Koch2019}
W.~Koch, R.~Mancuso, R.~West, and A.~Bestavros, ``Reinforcement learning for uav attitude control,'' {\em ACM Transactions on Cyber-Physical Systems}, vol.~3, no.~2, 2019.

\bibitem{Kaufmann2022}
E.~Kaufmann, L.~Bauersfeld, and D.~Scaramuzza, ``A benchmark comparison of learned control policies for agile quadrotor flight,'' in {\em ICRA 2022}, pp.~10504--10510, 2022.

\bibitem{Dhuheir2023}
M.~A. Dhuheir, E.~Baccour, A.~Erbad, S.~S. Al-Obaidi, and M.~Hamdi, ``Deep reinforcement learning for trajectory path planning and distributed inference in resource-constrained uav swarms,'' {\em IEEE Internet of Things Journal}, vol.~10, no.~9, pp.~8185--8201, 2023.

\bibitem{Azar2021}
A.~T. Azar, A.~Koubaa, N.~Ali~Mohamed, H.~A. Ibrahim, Z.~F. Ibrahim, M.~Kazim, A.~Ammar, B.~Benjdira, A.~M. Khamis, I.~A. Hameed, and G.~Casalino, ``Drone deep reinforcement learning: A review,'' {\em Electronics (Switzerland)}, vol.~10, no.~9, pp.~1--30, 2021.

\bibitem{Kurunathan2024}
H.~Kurunathan, H.~Huang, K.~Li, W.~Ni, and E.~Hossain, ``Machine learning-aided operations and communications of unmanned aerial vehicles: A contemporary survey,'' {\em IEEE Communications Surveys and Tutorials}, vol.~26, no.~1, pp.~496--533, 2024.

\bibitem{Bohn2019}
E.~Bohn, E.~M. Coates, S.~Moe, and T.~A. Johansen, ``Deep reinforcement learning attitude control of fixed-wing uavs using proximal policy optimization,'' {\em 2019 International Conference on Unmanned Aircraft Systems, ICUAS 2019}, pp.~523--533, 2019.

\bibitem{Bohn2024}
E.~Bohn, E.~M. Coates, D.~Reinhardt, and T.~A. Johansen, ``Data-efficient deep reinforcement learning for attitude control of fixed-wing uavs: Field experiments,'' {\em IEEE Transactions on Neural Networks and Learning Systems}, vol.~35, no.~3, pp.~3168--3180, 2024.

\bibitem{Xu2019}
J.~Xu, T.~Du, M.~Foshey, B.~Li, B.~Zhu, A.~Schulz, and W.~Matusik, ``Learning to fly: computational controller design for hybrid uavs with reinforcement learning,'' {\em ACM Trans. Graph.}, vol.~38, July 2019.

\bibitem{Muniraj2017}
D.~Muniraj, M.~C. Palframan, K.~T. Guthrie, and M.~Farhood, ``Path-following control of small fixed-wing unmanned aircraft systems with ${H}_\infty$ type performance,'' {\em Control Eng. Pract.}, vol.~67, pp.~76--91, 2017.

\bibitem{Yang2021}
J.~Yang, C.~Liu, M.~Coombes, Y.~Yan, and W.-H. Chen, ``Optimal path following for small fixed-wing uavs under wind disturbances,'' {\em IEEE Trans. Control Syst. Technol.}, vol.~29, no.~3, pp.~996--1008, 2021.

\bibitem{Bengio2009}
Y.~Bengio, J.~Louradour, R.~Collobert, and J.~Weston, ``Curriculum learning,'' in {\em Proc. of ICML}, pp.~41--48, 2009.

\bibitem{Tobin2017}
J.~Tobin, R.~Fong, A.~Ray, J.~Schneider, W.~Zaremba, and P.~Abbeel, ``Domain randomization for transferring deep neural networks from simulation to the real world,'' in {\em 2017 IEEE/RSJ International Conference on Intelligent Robots and Systems (IROS)}, pp.~23--30, 2017.

\bibitem{Hu2023}
B.~Hu, K.~Zhang, N.~Li, M.~Mesbahi, M.~Fazel, and T.~Başar, ``Toward a theoretical foundation of policy optimization for learning control policies,'' {\em Annual Review of Control, Robotics, and Autonomous Systems}, vol.~6, pp.~123--158, 2023.

\bibitem{Pattanaik2017}
A.~Pattanaik, Z.~Tang, S.~Liu, G.~Bommannan, and G.~Chowdhary, ``Robust deep reinforcement learning with adversarial attacks,'' {\em arXiv preprint arXiv:1712.03632}, 2017.

\bibitem{Huang2017}
S.~Huang, N.~Papernot, I.~Goodfellow, Y.~Duan, and P.~Abbeel, ``Adversarial attacks on neural network policies,'' {\em arXiv preprint arXiv:1702.02284}, 2017.

\bibitem{Muniraj2018}
D.~Muniraj, K.~G. Vamvoudakis, and M.~Farhood, ``Enforcing signal temporal logic specifications in multi-agent adversarial environments: A deep q-learning approach,'' {\em Proc. of the IEEE Conference on Decision and Control}, vol.~2018-December, pp.~4141--4146, 2018.

\bibitem{Lin2019}
Y.~Lin, Z.~Hong, Y.~Liao, M.~Shih, M.~Liu, and M.~Sun, ``Tactics of adversarial attack on deep reinforcement learning agents,'' 2019.

\bibitem{Pinto2017}
L.~Pinto, J.~Davidson, R.~Sukthankar, and A.~Gupta, ``Robust adversarial reinforcement learning,'' {\em 34th International Conference on Machine Learning, ICML 2017}, vol.~6, pp.~4310--4319, 2017.

\bibitem{Wang2024}
L.~Wang, S.~Zheng, S.~Tai, H.~Liu, and T.~Yue, ``Uav air combat autonomous trajectory planning method based on robust adversarial reinforcement learning,'' {\em Aerospace Science and Technology}, vol.~153, p.~109402, 2024.

\bibitem{Farhood2008}
M.~Farhood and G.~E. Dullerud, ``Control of systems with uncertain initial conditions,'' {\em IEEE Transactions on Automatic Control}, vol.~53, pp.~2646--2651, 12 2008.

\bibitem{schulman2017}
J.~Schulman, F.~Wolski, P.~Dhariwal, A.~Radford, and O.~Klimov, ``Proximal policy optimization algorithms,'' 2017.

\bibitem{Schulman2015}
J.~Schulman, S.~Levine, P.~Abbeel, M.~Jordan, and P.~Moritz, ``Trust region policy optimization,'' in {\em Proc. of ICML}, pp.~1889--1897, 2015.

\bibitem{Marquis2025}
D.~J. Marquis and M.~Farhood, ``Development and application of a dynamic obstacle avoidance algorithm for small fixed-wing aircraft with safety guarantees.'' Submitted.

\bibitem{Fry2020}
J.~M. Fry and M.~Farhood, ``A comprehensive analytical tool for control validation of fixed-wing unmanned aircraft,'' {\em IEEE Transactions on Control Systems Technology}, vol.~28, no.~5, pp.~1785--1801, 2020.

\bibitem{Real1993}
T.~R. Real, ``Digital simulation of atmospheric turbulence for dryden and von karman models,'' {\em Journal of Guidance, Control, and Dynamics}, vol.~16, no.~1, pp.~132--138, 1993.

\bibitem{Arifianto2015}
O.~Arifianto and M.~Farhood, ``Optimal control of a small fixed-wing {UAV} about concatenated trajectories,'' {\em Control Engineering Practice}, vol.~40, pp.~113--132, 2015.

\bibitem{Raffin2021StableBaselines3}
A.~Raffin, A.~Hill, A.~Gleave, A.~Kanervisto, M.~Ernestus, and N.~Dormann, ``Stable-baselines3: Reliable reinforcement learning implementations,'' {\em The Journal of Machine Learning Research}, vol.~22, no.~1, pp.~12348--12355, 2021.

\bibitem{paszke2019pytorch}
A.~Paszke, S.~Gross, F.~Massa, A.~Lerer, J.~Bradbury, G.~Chanan, T.~Killeen, Z.~Lin, N.~Gimelshein, L.~Antiga, {\em et~al.}, ``Pytorch: An imperative style, high-performance deep learning library,'' {\em Advances in Neural Information Processing Systems}, vol.~32, 2019.

\bibitem{ArifiantoFarhood_DSCC2012}
O.~Arifianto and M.~Farhood, ``Optimal control of fixed-wing {UAVs} along real-time trajectories,'' in {\em Proceedings of the ASME 2012 5th Annual Dynamic Systems and Control Conference}, pp.~205--214, 2012.

\bibitem{Lofberg2004}
J.~L{\"{o}}fberg, ``Yalmip : A toolbox for modeling and optimization in matlab,'' in {\em Proceedings of CACSD Conference}, (Taipei, Taiwan), 2004.

\bibitem{mosek}
M.~ApS, ``{The MOSEK optimization toolbox for MATLAB manual. Version 9.1.},'' 2019.

\end{thebibliography}

\appendix

\subsection{CZ-150 Model Details}
\label{app:model}

This section summarizes the CZ-150’s physical properties, aerodynamic model, and actuator model from \cite{Marquis2025}. The inertial and geometric properties of the CZ-150 are listed in Table~\ref{table:properties}. The aerodynamic force and moment coefficients, with values given in Table~\ref{table:coefficients}, are defined as
\begin{equation*}
\label{eq:aeromodel}
\begin{aligned}
C_X &= C_{X_{\alpha^2}} \alpha^2 + C_{X_{\mathcal{J}}} \mathcal{J} + C_{X_{\mathcal{J}^2}} \mathcal{J}^2 
     + C_{X_{\delta_E \alpha}} \delta_E \alpha + C_{X_0}, \\
C_Y &= C_{Y_{\beta}} \beta + C_{Y_{\hat{p}}} \hat{p} + C_{Y_{\hat{r}}} \hat{r} 
     + C_{Y_{\delta_A}} \delta_A + C_{Y_{\delta_R}} \delta_R + C_{Y_0}, \\
C_Z &= C_{Z_{\alpha}} \alpha + C_{Z_{\hat{q}}} \hat{q} 
     + C_{Z_{\delta_E}} \delta_E + C_{Z_0}, \\
C_L &= C_{L_{\beta}} \beta + C_{L_{\hat{p}}} \hat{p} + C_{L_{\hat{r}}} \hat{r} 
     + C_{L_{\delta_A}} \delta_A + C_{L_0}, \\
C_M &= C_{M_{\alpha}} \alpha + C_{M_{\alpha^3}} \alpha^3 + C_{M_{\hat{q}}} \hat{q} \\
    &\quad + C_{M_{\delta_E}} \delta_E + C_{M_{\hat{q}\delta_E}} \hat{q}\delta_E + C_{M_0}, \\
C_N &= C_{N_{\beta}} \beta + C_{N_{\hat{p}}} \hat{p} + C_{N_{\hat{r}}} \hat{r}\\
    &\quad + C_{N_{\delta_A}} \delta_A + C_{N_{\delta_R}} \delta_R + C_{N_0}.
\end{aligned}
\end{equation*}

\begin{table}[ht]
\centering
\caption{CZ-150 Inertial and Geometric Properties} 
\label{table:properties}
\begin{tabular}{lll|lll}
\toprule
\textbf{Property} & \textbf{Value} & \textbf{Units} & \textbf{Property} & \textbf{Value} & \textbf{Units} \\
\midrule
$m$ & 4.90 & kg & $J_{xx}$ & 0.546 & kg$\cdot$m$^2$ \\
$\bar{c}$ & 0.320 & m & $J_{yy}$ & 0.430 & kg$\cdot$m$^2$ \\
$b$ & 2.12 & m & $J_{zz}$ & 0.801 & kg$\cdot$m$^2$ \\
$S$ & 0.680 & m$^2$ & $J_{xz}$ & 0.066 & kg$\cdot$m$^2$ \\
$D_\mathrm{prop}$ & 0.406 & m & & & \\
\bottomrule
\end{tabular}
\end{table}

\begin{table}[ht]
 \centering
 \caption{CZ-150 Aerodynamic Coefficients}
 \label{table:coefficients}
 \begin{tabular}{lc|lc|lc}
 \toprule
 $\mathbf{C_X}$ & $\hat{\theta}$ & $\mathbf{C_Y}$ & $\hat{\theta}$ & $\mathbf{C_Z}$ & $\hat{\theta}$ \\
 \midrule
 $C_{X_{\alpha^2}}$ & +3.82 & $C_{Y_{\beta}}$ & -0.613 & $C_{Z_{\alpha}}$ & -4.925 \\ 
 $C_{X_{\mathcal{J}}}$ & +0.111 & $C_{Y_{\hat{p}}}$ & -0.136 & $C_{Z_{\hat{q}}}$ & +16.9 \\ 
 $C_{X_{\mathcal{J}^2}}$ & +0.0575 & $C_{Y_{\hat{r}}}$ & -0.284 & $C_{Z_{\delta_E}}$ & -0.161 \\ 
 $C_{X_{\delta_E \alpha}}$ & +1.08 & $C_{Y_{\delta_A}}$ & -0.131 & $C_{Z_{0}}$ & +0.0296 \\ 
 $C_{X_{0}}$ & -0.00680 & $C_{Y_{\delta_R}}$ & +0.0481 & \phantom{} & \phantom{} \\ 
 \phantom{} & \phantom{} & $C_{Y_0}$ & +0.0214 & \phantom{} & \phantom{} \\
 \midrule
 $\mathbf{C_L}$ & $\hat{\theta}$ & $\mathbf{C_M}$ & $\hat{\theta}$ & $\mathbf{C_N}$ & $\hat{\theta}$ \\
 \midrule
 $C_{L_{\beta}}$ & -0.0530 & $C_{M_{\alpha}}$ & -0.356 & $C_{N_{\beta}}$ & +0.0390 \\ 
 $C_{L_{\hat{p}}}$ & -0.215 & $C_{M_{\alpha^3}}$ & +1.85 & $C_{N_{\hat{p}}}$ & +0.00470 \\ 
 $C_{L_{\hat{r}}}$ & +0.0326 & $C_{M_{\hat{q}}}$ & -1.59 & $C_{N_{\hat{r}}}$ & -0.0991 \\ 
 $C_{L_{\delta_A}}$ & -0.0758 & $C_{M_{\delta_E}}$ & -0.197 & $C_{N_{\delta_A}}$ & +0.0150 \\ 
 $C_{L_{0}}$ & -0.0002 & $C_{M_{\hat{q}\delta_E}}$ & +9.71 & $C_{N_{\delta_R}}$ & -0.0259 \\ 
 \phantom{} & \phantom{} &$C_{M_0}$ & +0.0340 & $C_{N_0}$ & +0.00004 \\ 
 \bottomrule
 \end{tabular}
\end{table}

Commanded actuator deflections are derived from centered PWM signals using static polynomial mappings:
\begin{equation*}
\label{eq:actuator_static}
 \delta_i^\mathrm{cmd} =
 \begin{cases} 
 \dfrac{\pi}{180}\displaystyle\sum_{k=0}^{3} c_{k,i} (\delta_i^\mathrm{PWM})^k, & i \in \{E, A, R\}, \\[2ex]
 \dfrac{1}{60}\left(c_{1,i}\delta_i^\mathrm{PWM} + c_{0,i}\right), & i = T,
 \end{cases}
\end{equation*}
where $c_{k,i}$ denotes the coefficient of the $k^\text{th}$ polynomial term and the coefficients are given in Table~\ref{table:actuator_dynamic}. Each actuator follows a first-order dynamic response:
\begin{equation*}
\dot{\delta}_i = \frac{1}{\tau_i}(\delta_i^\mathrm{cmd} - \delta_i), \quad i \in \{E, A, R, T\},
\end{equation*}
where the time constants are $\tau_A = 0.083$\,s, $\tau_E = 0.071$\,s, $\tau_R = 0.071$\,s, and $\tau_T = 0.082$\,s, respectively.

\begin{table}[ht]
\centering
\caption{CZ-150 Coefficients for Control Surfaces}
\label{table:actuator_dynamic}
\begin{tabular}{c c| c c}
\toprule
\textbf{Coefficient} & \textbf{Value} & \textbf{Coefficient} & \textbf{Value} \\
\midrule
 $c_{3,A}$ & $-3.98\times10^{-8}$ & $c_{2,A}$ & $1.08\times10^{-8}$ \\
 $c_{1,A}$ & $-6.24\times10^{-2}$ & $c_{0,A}$ & $-2.56\times10^{-1}$ \\
\midrule
 $c_{3,E}$ & $-8.24\times10^{-8}$ & $c_{2,E}$ & $-9.54\times10^{-6}$ \\
 $c_{1,E}$ & $7.21\times10^{-2}$ & $c_{0,E}$ & $-8.87\times10^{-1}$ \\
\midrule

 $c_{3,R}$ & $-2.72\times10^{-7}$ & $c_{2,R}$ & $7.58\times10^{-6}$ \\
 $c_{1,R}$ & $1.16\times10^{-1}$ & $c_{0,R}$ & $-7.40\times10^{-1}$ \\
\midrule
$c_{1,T}$ & $6.56\times10^3$ & $c_{0,T}$ & $2.91\times10^3$ \\
\bottomrule
\end{tabular}
\end{table}

The following aerodynamic coefficient perturbation bounds were obtained from system identification:
\begin{equation*}
\begin{aligned}
\begin{bmatrix} \boldsymbol{\Delta}_{\mathbf{C}}^- & \boldsymbol{\Delta}_{\mathbf{C}}^+ \end{bmatrix} &=
\begin{bmatrix}
-0.0258 & 0.0258\\ -0.0510 & 0.0510\\ -0.0872 & 0.0872 \\ -0.0204 &0.0204 \\ -0.0330 & 0.0330 \\ -0.0084 & 0.0084
\end{bmatrix}, \\
\begin{bmatrix} \boldsymbol{\nu}_{\mathbf{C}}^- & \boldsymbol{\nu}_{\mathbf{C}}^+ \end{bmatrix} &=
\begin{bmatrix}
-0.0180 & 0.0175 \\
-0.0289 & 0.0287 \\
-0.0598 & 0.0606 \\
-0.0128 & 0.0128 \\
-0.0299 & 0.0299 \\
-0.0044 & 0.0044
\end{bmatrix}.
\end{aligned}
\end{equation*}

\subsection{Switched UIC Controller Design}
\label{app:uic}
The switched UIC controller used in this work is based on the control synthesis technique developed in~\cite{Farhood2008} for discrete-time LTI systems with uncertain initial conditions. For trajectories composed of concatenated trim segments, a switched plant formulation can capture the changing dynamics, and the UIC approach further accounts for uncertain initial states arising from switching~\cite{ArifiantoFarhood_DSCC2012,Arifianto2015}. A UIC subcontroller is typically $N$-eventually time invariant, meaning that its state-space matrix sequences are all $N$-eventually time invariant for some integer $N \geq 1$, e.g., $(A_0^K, A_1^K, \dots, A_{N-1}^K, A_N^K, A_N^K, \dots)$. Given a switched system representing the plant, the control objective is to design an $N_i$-eventually time-invariant subcontroller $K_i$ for each constituent model $P_i$ of the switched plant, for some integer $N_i\geq 1$, such that the closed-loop system composed of  $P_i$ and $K_i$ is asymptotically stable and the performance output $\mathbf{z}$ satisfies
\begin{equation*}
\sup\{||\mathbf{z}||_{\ell_2} \mid ||\tilde{\mathbf{x}}|| \le 1, \; ||\mathbf{d}||_{\ell_2} \le 1\} < \xi_\mathrm{min},
\end{equation*}
where $\tilde{\mathbf{x}}$ is a vector of uncertain initial state values, $\mathbf{d}$ is the disturbance input, and $\xi_\mathrm{min}$ denotes the minimum achievable performance level up to a specified tolerance. Here, $\ell_2$ is the Hilbert space of sequences $\mathbf{v}=(\mathbf{v}_0,\mathbf{v}_1,\ldots)$ with $\mathbf{v}_k \in \mathbb{R}^{n_v}$ such that $\sum_{k=0}^{\infty} \mathbf{v}_k^T \mathbf{v}_k < \infty$, and the $\ell_2$-norm of $\mathbf{v}$ is defined as $||\mathbf{v}||_{\ell_2}^2 = \sum_{k=0}^{\infty} \mathbf{v}_k^T \mathbf{v}_k$.

The nonlinear equations of motion of the aircraft are linearized about various trim trajectories corresponding to nominal airspeed $V_{\mathrm{r,nom}} = 21$\,m/s and $(\kappa,\gamma) \in K_\mathrm{con} \times \Gamma_\mathrm{con}$, where $K_\mathrm{con} = \{-0.02, -0.012, 0, 0.012, 0.02\}$ and $\Gamma_\mathrm{con} = \{-0.21, -0.11, 0, 0.11, 0.21\}$, yielding 25 linearized models. Trim values are obtained as described in Section~\ref{sec:ref}, and each linearized model is discretized using zero-order hold sampling with a sampling period of $\Delta t = 0.04$\,s. A finite horizon of $N=5$ is used for all subcontrollers.

The performance output $\mathbf{z}$ is defined in terms of the tracking errors in attitude, angular rates, position, and actuator deflections:
\begin{equation*}
\begin{aligned}
\mathbf{z} =
\biggl[
\frac{1}{2\pi}\bar{\boldsymbol{\omega}}^T \ \ 0.1 \bar{V}_r \ \ \frac{4}{\pi} \bar{\phi} \ \ \frac{4}{\pi} \bar{\theta} \ \ \frac{1}{\pi} \bar{\psi} \ \ 0.2 \bar{\mathbf{p}}^T \ \ 5 \bar{\boldsymbol{\delta}}^T\biggr]^T.
\end{aligned}
\end{equation*}
For concatenated trajectories, each segment corresponds to a specific $(\kappa,\gamma)$ pair, so transitioning to a new segment induces a switch in the plant dynamics. Nonlinear trim analysis indicates that variations in $\kappa$ primarily affect yaw rate ($r$) and bank angle ($\phi$), while variations in $\gamma$ predominantly influence pitch angle ($\theta$), elevator deflection ($\delta_E$), and throttle ($\delta_T$). These relationships inform the selection of the matrix $\Lambda \in \mathbb{R}^{16 \times 5}$ used to express the nonzero uncertain initial state $\bar{\mathbf{x}}_0 = \Lambda \tilde{\mathbf{x}}$. The nonzero entries of $\Lambda$ are $\Lambda_{3,1} = 0.1$, $\Lambda_{7,2} = 0.5$, $\Lambda_{8,3} = 0.5$, $\Lambda_{13,4} = 0.2$, and $\Lambda_{16,5} = 0.2$. Controller synthesis proceeds in two stages. First, a semidefinite program (SDP) is solved to determine the minimum achievable performance level $\xi_\mathrm{min}$, following the formulation in~\cite{Farhood2008}. Then, a second SDP is solved to generate a realizable controller that achieves a slightly relaxed performance level $\xi$ (here, $\xi = 1.3\,\xi_\mathrm{min}$) to ensure a well-conditioned solution. The SDPs are solved in MATLAB using the YALMIP modeling language~\cite{Lofberg2004} and the \texttt{MOSEK} solver~\cite{mosek}.
\end{document}